\documentclass[conference, a4paper]{IEEEtran}
\IEEEoverridecommandlockouts
\usepackage[a4paper, left=0.75in, right=0.51in, bottom=0.75in, top=1.69in]{geometry}
\usepackage{cite}
\usepackage{amsmath,amssymb,amsfonts}
\usepackage{graphicx}
\usepackage{textcomp}
\usepackage{xcolor}
\def\BibTeX{{\rm B\kern-.05em{\sc i\kern-.025em b}\kern-.08em
    T\kern-.1667em\lower.7ex\hbox{E}\kern-.125emX}}
\usepackage{tikz}
\usepackage{xcolor}
\usetikzlibrary{backgrounds}
\usepackage{setspace}
\usetikzlibrary{positioning}
\usetikzlibrary{shapes.geometric}
\usepackage{ifpdf}
\IEEEoverridecommandlockouts  
\usepackage[utf8]{inputenc}
\usepackage{amsmath}
\usepackage{graphicx}
\usepackage{flushend}
\usepackage{caption}
\usepackage{subcaption}
\usepackage{comment}
\usepackage{bm}
\usepackage{xspace}
\usepackage{soul}
\usepackage{amsopn}
\usepackage{algpseudocode}
\usepackage{algorithm} 
\usepackage{enumerate}
\makeatletter
\let\NAT@parse\undefined
\makeatother
\usepackage{amsmath,amsfonts, amssymb,amsthm,color}
\usepackage[numbers,sort&compress]{natbib}
\usepackage{amsf onts}
\usepackage{color}
\usepackage{setspace}

\usepackage{amssymb}
\newtheorem{theorem}{Theorem}

\usepackage{breqn}
\usepackage{hyperref} 
\usepackage{epstopdf}
\usepackage{lipsum}
\usepackage{epsfig}
\DeclareGraphicsExtensions{.pdf,.eps,.png,.jpg,.mps}

\usepackage{hyperref}
\hypersetup{
    colorlinks=true,
    linkcolor=blue,
    filecolor=magenta,      
    urlcolor=cyan,
}
\usepackage{rotating}

\begin{document}
\title{
Comparison of Minimization Methods for  Rosenbrock Functions
}
\author{\IEEEauthorblockN{Iyanuoluwa Emiola}
\IEEEauthorblockA{\textit{Electrical and Computer Engineering} \\
\textit{University of Central Florida}\\
Orlando FL 32816, USA \\
iemiola@knights.ucf.edu}
\and
\IEEEauthorblockN{Robson Adem}
\IEEEauthorblockA{\textit{Electrical and Computer Engineering} \\
\textit{University of Central Florida}\\
Orlando FL 32816, USA \\
ademr@knights.ucf.edu}
}


\maketitle


\begin{abstract}
This paper gives an in-depth review of the most common iterative methods for unconstrained optimization using two functions that belong to a class of Rosenbrock functions as a performance test. This study covers the Steepest Gradient Descent Method, the Newton-Raphson Method, and the Fletcher-Reeves
Conjugate Gradient method. In addition, four different step-size selecting methods including fixed-step-size, variable step-size, quadratic-fit, and golden section method were considered. Due to the computational nature of solving minimization problems, testing the algorithms is an essential part of this paper. Therefore, an extensive set of numerical test results is also provided to present an insightful and a comprehensive comparison of the reviewed algorithms. This study highlights the differences and the trade-offs involved in comparing these algorithms.
 \end{abstract}
 \begin{IEEEkeywords}
Rosenbrock functions, gradient descent methods, variable step-size, quadratic fit, golden section method, conjugate gradient method, Newton-Raphson.
\end{IEEEkeywords}
 \section{Introduction}
 Solutions to unconstrained optimization problems can be applied to multi-agent systems and  machine learning problems, especially if the problem is in a decentralized or distributed fashion \cite{nedic2018network}, \cite{yang2016distributed}, \cite{montijano2014efficient} and \cite{tsianos2012consensus}. Sometimes, in adversarial attack applications, malicious agents can be present in a network that will slow down convergence rates to optimal points as seen in \cite{emiola2021distributed}, \cite{sundaram2018distributed} and \cite{ravi2019case}. Therefore the need for a fast convergence and cost associated with it are usually  necessities by recent researchers. The steepest descent method is a good first order method for obtaining optimal solutions if an appropriate step size is chosen. Some methods of choosing step sizes include the fixed step size, the variable step size, the polynomial fit, and the golden section method that will be discussed in details in subsequent sections. Nonetheless, these methods have their own merits and demerits. A second order method such as the Newton-Raphson method is very suitable for quadratic problems and attains optimality in a small number of iterations \cite{chong2004introduction}. However the Newton-Raphson method requires computing the Hessian and its inverse which are often a bottleneck. For this reason, the Newton-Raphson  method may not be suitable for solving large scale optimization problems. 
 
 To address the lapses that the Newton-Raphson method poses, some methods have been recently proposed such as matrix splitting discussed in \cite{zargham2013accelerated}. Another way of obtaining fast convergence properties of second-order methods using the structure of first-order methods are the Quasi-Newton methods.  These methods  incorporate second-order (curvature information) in the first-order approaches. Examples of these methods include the Broyden-Fletcher-Goldfarb-Shanno algorithm (BFGS) \cite{eisen2017decentralized} and the Barzilai-Borwein (BB) \cite{dai2005projected} and \cite{gill1972quasi}. However these methods usually require additional assumptions on the objective function to be minimized to be strongly convex and that the gradient of such function  to be Lipschitz continuous to improve convergence rate discussed in \cite{gao2019geometric}.
 The conjugate gradient method \cite{hestenes1952methods} on the other hand does not have as much restrictions as the Newton-Raphson method in terms of computing the inverse of the Hessian. The method computes the direction of search at each iteration and the direction is expressed as a linear combination of previous search direction calculation and the present gradient at the present iteration. 
The other interesting attribute of the conjugate gradient method is that it gives different methods of calculating the search directions and it is not only limited to quadratic functions as we will use that while performing simulations using the Rosenbrock function \cite{chong2004introduction}.  

The crux of our work entails comparing the convergence properties of first order and second order methods and the limitations of each method using two Rosenbrock functions. One of the two functions is usually known as the \textit{banana function} and is much more difficult to minimize than the other. We compare convergence attributes for these methods by examining these Rosenbrock functions where one is a alteration of the other. The first order type methods we used for this study include the steepest descent with fixed step size, variable step size, polynomial quadratic fit and golden section method. For the second order method, we consider the Newton-Raphson method and compare its performance with the Conjugate Gradient Methods. The above methods highlighted have their own advantages and disadvantages in terms of general performance, convergence, precision, convergence rate, and robustness. Depending upon the nature of the problem, performance design specifications, and available resources, one can select the most appropriate optimization method. In order to assist this effort, this study will highlight the differences and the trade-offs involved in comparing these algorithms. 
\subsection{Contributions} 
In this paper, we compare different types of methods as described to analyse convergence without restricting the analysis to quadratic functions. As we will see in subsequent sections, we examined first order methods using four different cases like the fixed step size, variable step size, quadratic fit and the golden section methods. In our analysis, we did not depend on the steepest descent with fixed step size commonly used by many researchers to study first order methods. This paper also compares second order methods and highlights the advantages the Newton-Raphson and conjugate gradient methods have over each other. 
  \subsection{Paper Pattern}
Section \ref{sec:problemformulation} presents the problem formulation, section introduces the different types of first and second type methods, and their convergence analysis. Numerical experiments are performed in section \ref{sec:Numerical} and the conclusion is given in section \ref{sec:conclusion}. 
\subsection{Notation}
We denote the set of positive and negative reals as $\mathbb{R}_+$ and $\mathbb{R}_-$, the transpose of a vector or matrix as $(\cdot)^T$, and the L$2$-norm of a vector by $||\cdot||$. We let the gradient of a function $f(\cdot)$ to be $\nabla f(\cdot)$, the Hessian of a function $f(\cdot)$ be $F(\cdot)= \nabla^{2} f(\cdot)$, and $\langle\ \cdot ,\cdot \rangle$ denotes the inner product of two vectors. 
\section{Problem formulation}\label{sec:problemformulation}
In mathematical optimization, Rosenbrock functions are used as a performance test problem for optimization algorithms \cite{Rosenbrock}. So we will use the following  Rosenbrock functions in our analysis.

\begin{align}\label{Rosenbrock2}
&\underset{x_1,x_2}{\operatorname{minimize}}& & f(x_1,x_2)= \kappa(x_1^{2}-x_2)^{2}+(x_1 - 1)^{2}  
\end{align}
where $\kappa>0$, $x$ is a vector such that $x=[x_1, x_2]^{T}$ and $f(x)$ in problem \eqref{Rosenbrock2} is strictly convex and twice differentiable..
We will examine the key role $\kappa$ plays in affecting convergence in section \ref{sec:Numerical}. To solve  the minimization problems \eqref{Rosenbrock2} using gradient methods, we let the iterative equation generated by the minimization problem \eqref{Rosenbrock2} be given by:
\begin{equation}\label{regularagentminimization}
     x(k+1) = x(k) - \alpha(k) \nabla f(x(k)),
\end{equation}
where $k$ is the iteration, $\alpha(k)>0$ is the step size and $\nabla f(x(k))$ is the gradients of $f$ at each iterate $x(k)$. 
Different first order and second order methods for solving problem \eqref{Rosenbrock2} will be explored in the later part of the paper. 

\begin{figure}[ht]
    \centering
  \includegraphics[scale=0.52]{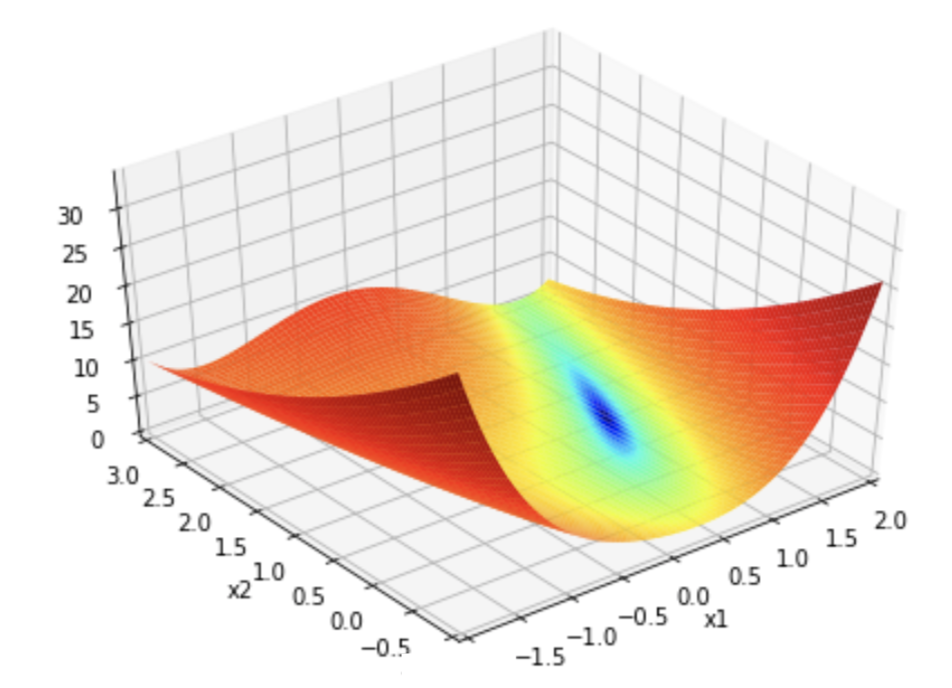}
  \caption{Surface plot for function (\ref{Rosenbrock2}) over $\mathbb {R} ^{2}$ when $\kappa=1$}
  \label{surface1}
\end{figure}

\begin{figure}[ht]
    \centering
  \includegraphics[scale=0.52]{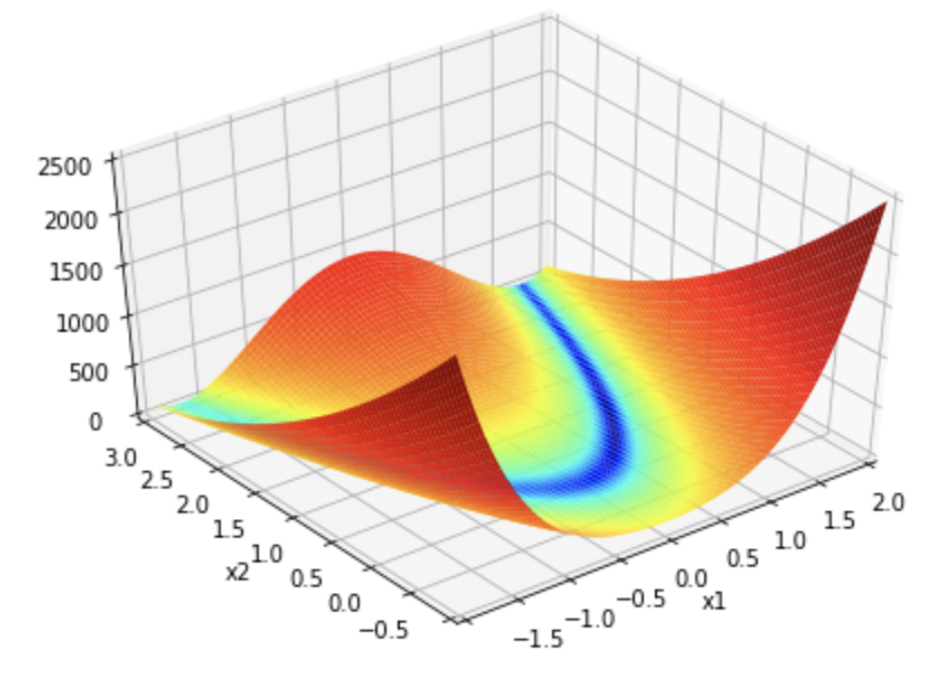}
  \caption{Surface plot for function (\ref{Rosenbrock2}) over $\mathbb {R} ^{2}$ when $\kappa=100$}
  \label{surface2}
\end{figure}
\begin{figure}[ht]
    \centering
  \includegraphics[scale=0.60]{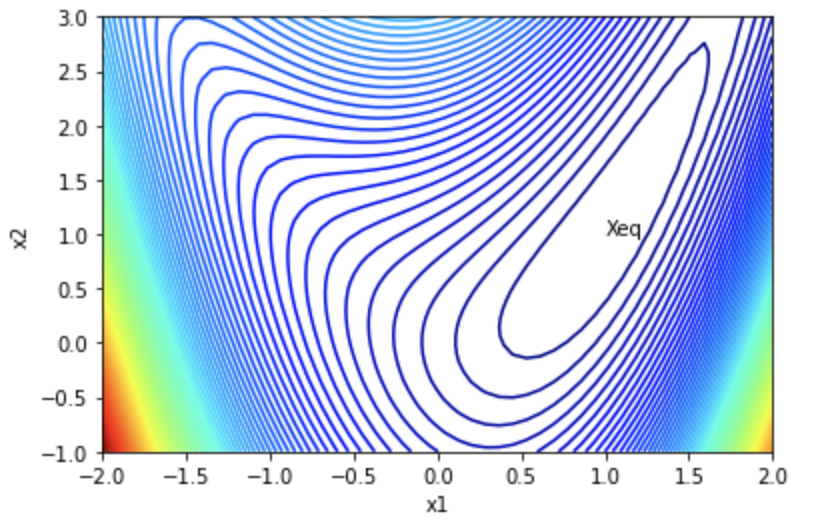}
  \caption{Level  curves for function (\ref{Rosenbrock2}) over $\mathbb {R} ^{2}$ when $\kappa=1$}
  \label{contour1}
\end{figure}

\begin{figure}[ht]
    \centering
  \includegraphics[scale=0.60]{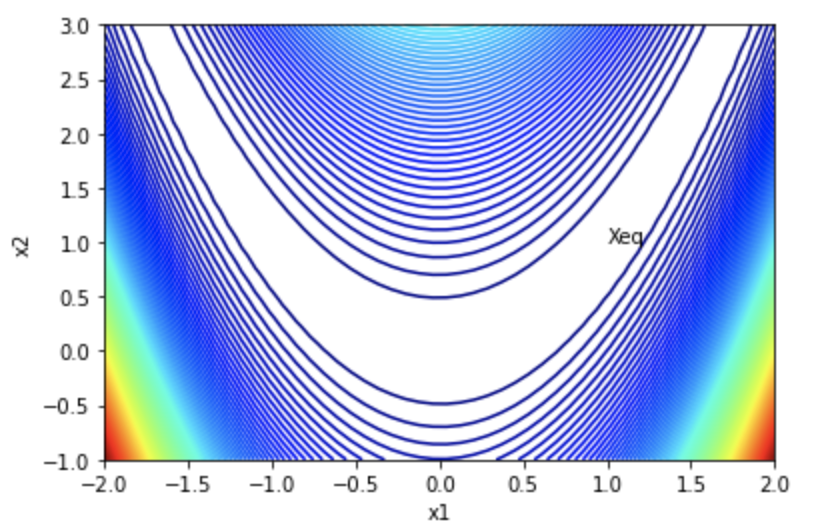}
  \caption{Level  curves for function (\ref{Rosenbrock2}) over $\mathbb {R} ^{2}$ when $\kappa=100$}
  \label{contour2}
\end{figure}
The level curves for function \ref{Rosenbrock2} with $\kappa=1$ and $\kappa=100$ are shown in Fig. \ref{contour1} and  Fig. \ref{contour2} respectively. Even though $f(x)$ has a global minimum at $(1,1)$, the global minimum of $f(x)$ when $\kappa=100$ is inside a longer blue parabolic shaped flat valley than it is when $\kappa=1$ as shown in Fig. \ref{surface1} and Fig. \ref{surface2}. Due to the longer size of the valley in Fig. \ref{surface2}, attaining convergence to the minimum of $f(x)$ with $\kappa=100$ becomes increasingly difficult. As such, the functions generated by \eqref{Rosenbrock2} with $\kappa=1$ and $\kappa=100$ are excellent candidates to evaluate the characteristics of optimization algorithms, such as: convergence rate, precision, robustness, and general performance. In this paper, function (\ref{Rosenbrock2}) with $\kappa=1$ and $\kappa=100$ will be  employed to establish a numerical comparison of the three optimization algorithms discussed below in section \ref{Firstandsecondorder}.

\section{Analysis of First and Second order methods}\label{Firstandsecondorder}
We now discuss the most common iterative methods for unconstrained optimization.
\subsection{The Steepest Descent Method}
To solve problem \eqref{Rosenbrock2}, a sequence of guesses $x(0), x(1),...x(k), x(k+1)$ will be generated in a descent manner such that $f(x(0))>f(x(1))>...>f(x(k+1))$. It can be often tedious to obtain optimality after some $K$ iterations, and $K$ is the maximum number of iterations needed for convergence such that $\nabla f(x(K))=0$. Therefore it suffices to actually modify the gradient stopping condition to satisfy $\|\nabla f(x(K))\|\leq\varepsilon$ where $\varepsilon>0$ and very small which is often referred to as the stopping criterion for convergence to hold. Different ways of choosing the step size will be explored below:
\subsubsection{Steepest Descent With a Constant Step Size}
The constant step size is constructed in a manner where you simply use one value of $\alpha$ in all iterations. 
To illustrate the fixed step size principle in solving problem \eqref{Rosenbrock2}, we will pick an  $\alpha$ value between $0$ and $1$ and show numerically how convergence is attained. 
Other methods of choosing the step size in a steepest descent algorithm are discussed below:
\subsubsection{Steepest Descent with Variable Step Size}
In the variable step size method, $3$ or $4$ values of $\alpha$ are chosen at each iteration and the value that produces the smallest $g(\alpha(k))$ value will be chosen where $g(\alpha(k))=f(x(k+1))$.
The variable step size algorithm is also easy to implement and has a better convergence probability than the fixed step size method. The results of simulating problem \eqref{Rosenbrock2} with $\kappa=1$ and $\kappa=100$ using the variable step size are shown in \eqref{sec:Numerical}.
\subsubsection{Steepest Descent with Quadratic Fit Method}
For the quadratic fit method, three values of $\alpha(k)$ are guessed at each iteration and the values of the corresponding $g(\alpha(k))$ values are computed, where $g(\alpha(k))=f(x(k+1))$ For example, suppose the three values of $\alpha$ values chosen are $\alpha(1), \alpha(2), \alpha(3)$. To fit a quadratic model of the form:
\begin{equation}\label{generalquadratic}
    g(\alpha)=a\alpha^2 + b\alpha +c,
\end{equation}
we write the quadratic model based on $\alpha$ values as:
\begin{equation}\label{firstquadratic}
    g(\alpha(1))=a\alpha(1)^2 + b\alpha(1) +c,
\end{equation}
\begin{equation}\label{secondquadratic}
    g(\alpha(2))=a\alpha(2)^2 + b\alpha(2) +c, 
\end{equation}
\begin{equation}\label{thirdquadratic}
      g(\alpha(3))=a\alpha(3)^2 + b\alpha(3) +c.
\end{equation}
where $a,b,c$ are constants. After solving for $a,b,c$ in equations \eqref{firstquadratic}, \eqref{secondquadratic}, and \eqref{thirdquadratic}, we will use these values in equation \eqref{generalquadratic} to obtain the optimum step size which is the mimimum of  equation \eqref{generalquadratic}.
  \subsubsection{Steepest Descent with the Golden Section Search}
  In this algorithm, we use a range between two values and divide the range into sections. We then eliminate some sections within the sections in the range to shrink the region where the convergence might occur. For this algorithm to be implemented as we will see in section \eqref{sec:Numerical}, the initial region of uncertainty and the stopping criterion have to be defined. An example where a golden section search was applied to minimize a function in a closed interval is seen in \cite{chong2004introduction}.
  
  Second-order methods have been an improvement in terms of convergence speeds when it comes to solving unconstrained optimization problems such as problem \eqref{Rosenbrock2}. We will show by simulations in section \ref{sec:Numerical} the speed of convergence for the Newton-Raphson and conjugate gradient methods compared to other methods. We will now analyze the two second order methods below:
\subsection{Newton-Raphson Methods}
The Newton-Raphson Method is very useful in obtaining fast convergence of an unconstrained problem like in equation \eqref{Rosenbrock2} especially when the initial starting point is very close to the minimum. The main disadvantage of this method is the cost and difficulty associated with finding the inverse of the Hessian and also ensuring that the Hessian inverse matrix is positive definite. The update equation for the Newton-Raphson Method is given by:
\begin{equation}\label{eqn:newton}
x(k+1) = x(k) - F^{-1}(x(k)){\nabla f(x(k))}.
\end{equation}
Some of the methods of approximating the term that contains the inverse of the Hessian, $ F^{-1}(x(k))$ in equation \eqref{eqn:newton} are the Quasi-Newton methods such as the BFGS and the Barzilai-Borwein methods \cite{eisen2017decentralized}, \cite{emiola2021sublinear}.
\subsection{Conjugate Gradient Methods}
For the class of quadratic functions $f(x) = 0.5x^{T}Qx-x^{T}b$,
and $x\in \mathbb{R}^{n}$,
the conjugate gradient algorithm uses a direction expressed in terms of the current gradient and the previous direction at each iteration by ensuring that the directions are mutually Q-conjugate, where $Q$ is a positive definite symmetric $n\times n$ matrix \cite{hestenes1952methods}. We note that the directions $d(0),d(1),.....d(m)$ are Q-conjugate if $d(i)^{T}Qd(j)=0$, and $i\neq j$. The conjugate gradient method also exhibits fast convergence property for non-quadratic problems like problem \eqref{Rosenbrock2}. In the simulation in section \ref{sec:Numerical}, we use the Fletcher-Reeves Formula \cite{fletcher1964function} given by:
\begin{equation*}
  \beta(k) =  \frac{g(k+1)^{T}g(k+1)}{g(k)^{T}g(k)}
\end{equation*}
where $g(k)=\nabla f(x(k))$ and $\beta(k)$ are constants picked such that the directional iteration $d(k+1)$ is Q-conjugate to $d(0),d(1),.....d(k)$ according to the following iterations:
\begin{equation*}
   x(k+1) = x(k) + \alpha(k)d(k),
\end{equation*}
and $\alpha(k)>0$ is the step size.
We will show in section \ref{sec:Numerical} that the conjugate gradient method performs better than the steepest descent method in terms of convergence rate when we use the same fixed step size.

\section{Numerical Experiments Insights}\label{sec:Numerical}
\begin{figure*}
  \includegraphics[width=\textwidth]{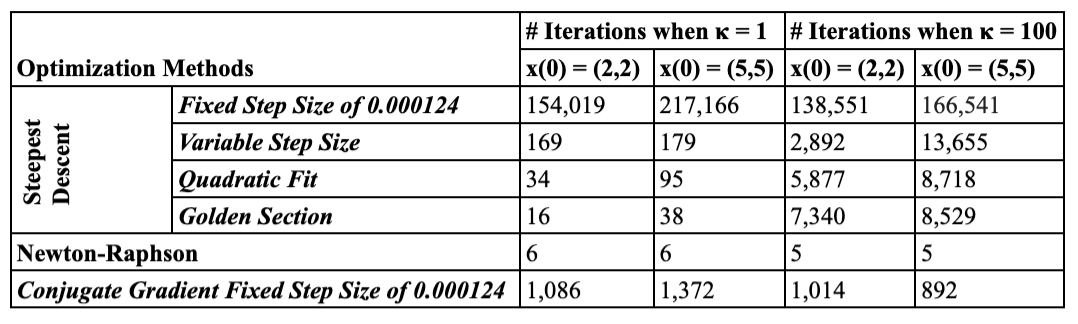}
  \caption{Numerical Comparison of the three optimization methods for the two functions with step-size = 0.000124 }
  \label{table}
\end{figure*}

\begin{figure}[ht]
  \includegraphics[scale = 0.345]{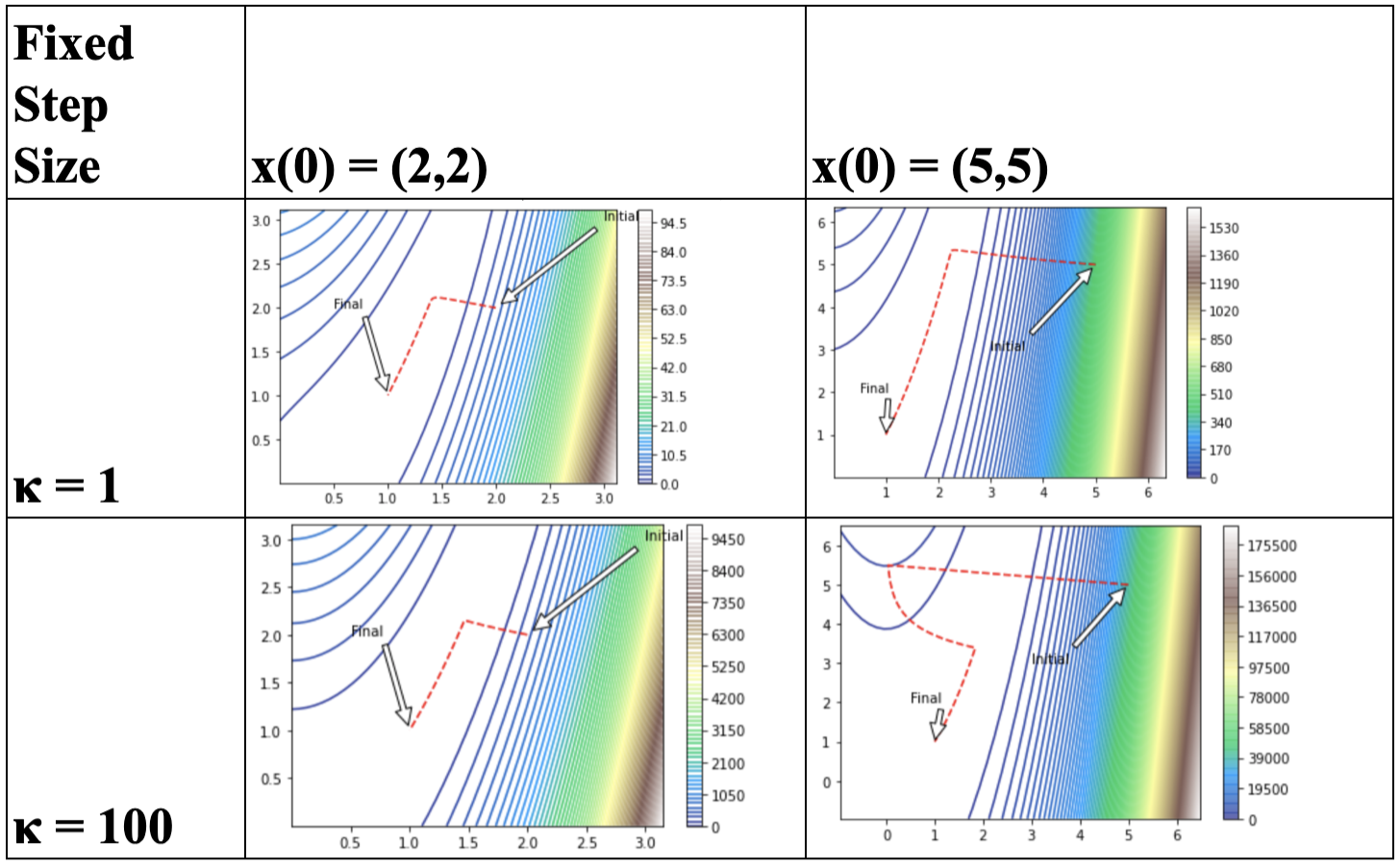}
  \caption{Level curves for Steepest Gradient Descent with fixed step size.}
  \label{fixed}
\end{figure}

\begin{figure}[ht]
  \includegraphics[scale = 0.247]{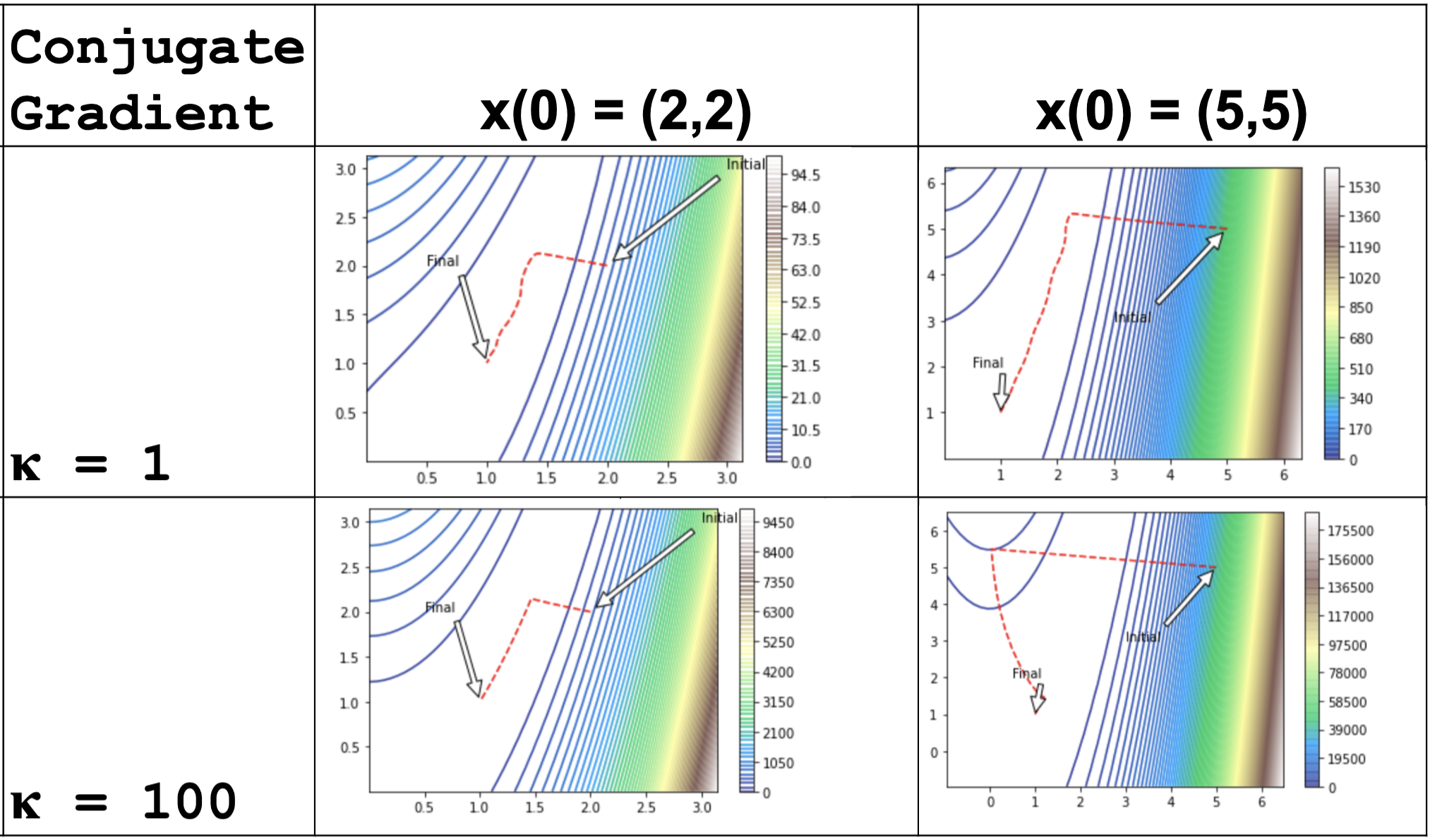}
  \caption{Level curves for Conjugate Gradient method.}
  \label{CG}
\end{figure}

In this section, we will compare the methods discussed in section \ref{Firstandsecondorder} in terms of convergence and the number of iterations taken to reach optimality in problem \eqref{Rosenbrock2} for cases when $\kappa =1$ and $\kappa =100$.  For these two cases, initial conditions of $(2,2)$ and $(5,5)$ and a stopping criterion of $\|\nabla f(x(K))\| \leq 0.001$ are used across all these methods discussed where $K$ is the maximum iteration number to achieve convergence. Fig. \ref{table} summarizes the numerical result whereas Fig. \ref{fixed} - Fig \ref{CG} illustrate the trajectory of iterates $x(k)$ on level curves of the functions across all cases. The code for the experiments is available \href{https://github.com/robsonadem/Comparison-of-Minimization-Methods-for-Rosenbrock-Functions}{here.}


\subsection{Geometric Step Size Test}
To compare the methods discussed in section \ref{Firstandsecondorder}, a geometric sequence of fixed step sizes is used.  We note that the Newton-Raphson method is not dependent on any step-size. As a result, the geometric step size test includes the the steepest descent and Conjugate Gradient Methods. Using the fixed step sizes, $0.124$, $0.0124$, $0.00124$ and $0.000124$, we present the results as follows. 

\subsubsection{Case when \texorpdfstring{$\alpha=0.124$}{Lg}}
When $\kappa=1$ and the step size $\alpha=0.124$ is used, the minimization problem \eqref{Rosenbrock2} generated by \eqref{regularagentminimization}  converges for the steepest descent with the $(2,2)$ initial condition and diverges with the initial condition of $(5,5)$. Equation \eqref{regularagentminimization} also diverges with the conjugate gradient with both of the initial starting points.

When $\kappa=100$ and the step size $\alpha =0.124$ is used, the problem \eqref{Rosenbrock2} generated by \eqref{regularagentminimization} diverges by the steepest descent for both initial points, and also diverges by the conjugate gradient using both initial points, $(2,2)$ and $(5,5)$.
\subsubsection{Case when \texorpdfstring{$\alpha=0.0124$}{Lg}}
When $\kappa=1$, and the steepest descent method is applied on the function \eqref{Rosenbrock2} generated by iteration \eqref{regularagentminimization}, convergence is attained with both of the starting initial points. This is an improvement over the case when $\alpha=0.124$, where the steepest method diverges using the initial point $(5,5)$. By the conjugate method on  function, convergence is obtained with both starting points which is also an improvement over the case when $\alpha=0.124$, where it diverges for both starting values.

For equation \eqref{Rosenbrock2} generated by \eqref{regularagentminimization} when the step size $\alpha =0.0124$ is used with $\kappa=100$, there was no improvement in convergence attributes because divergence is obtained by steepest descent and conjugate gradient methods for both of the initial points.

\subsubsection{Case when \texorpdfstring{$\alpha=0.00124$}{Lg}}
At the stage when $\alpha = 0.00124$ is used, the conjugate-gradient and the steepest descent both converge with both starting points for function \eqref{Rosenbrock2} generated by \eqref{regularagentminimization} with $\kappa=1$. 
By using equation \eqref{Rosenbrock2} generated by \eqref{regularagentminimization} with $\kappa=100$ and $\alpha =0.00124$, convergence is obtained using the initial point $(2,2)$ compared to divergence result obtained by using the initial point $(5, 5)$.

\subsubsection{Case when \texorpdfstring{$\alpha=0.000124$}{Lg}}
When the step size $\alpha=0.000124$ is used, the steepest descent method and conjugate gradient method converge for  equation \eqref{Rosenbrock2} generated by \eqref{regularagentminimization} for both $\kappa=1$ and $\kappa=100$. In addition, for each function, we observed that for both initial points $(2,2)$ and $(5,5)$ result in convergence. Therefore we will use this step
size as a case study to compare the rate of convergence of steepest descent and conjugate gradient for the two functions with the two starting points. 

\subsection{Significance of the Newton-Raphson Method}
The significance of the Newton-Raphson method should not be overlooked even for non-quadratic functions like \eqref{Rosenbrock2} because the method guarantees convergence to the optimal solution without specifying a step size. Moreover, it achieves convergence in just few iterations for both of the starting points as well as for the two functions. This affirms the unique convergence attribute of the Newton-Raphson method when the starting point is not far away from the optimal solution. 


\subsection{Comparison of the Variable Step Size, Quadratic Fit and Golden Section with other Methods} 
Starting with the two initial points $(2,2)$ and $(5,5)$ using the variable step size method, convergence was obtained for the function \eqref{Rosenbrock2} for both cases $\kappa=1$ and $\kappa=100$ when three varying step sizes of $0.000124$, $0.0124$ and $0.124$ are used. When the quadratic fit is used, three values of the step sizes are used in each iteration and selected from the range $(0.00001, 0.000124)$. The result from the quadratic fit shows that a better convergence is achieved for function \eqref{Rosenbrock2} with $\kappa=1$ but shows a weaker convergence for the second function. This explains that fluctuations in the random selection of step sizes between a range can influence the convergence rate. Moreover, the alteration parameter $\kappa$ from function \eqref{Rosenbrock2} can also slow down convergence rates. For the golden section method, a range of $(0.00000124,1.5)$ is used to locate the value of the step size that result in the solution to the minimization problem \eqref{Rosenbrock2}. By using the initial points of $(2,2)$ and $(5,5)$ on the function \eqref{Rosenbrock2} with $\kappa=1$, the golden section method resulted in the fastest convergence rate by comparing with the steepest descent with fixed step size, variable step size and the quadratic fit methods. However, convergence with the golden section is faster than the variable, fixed step size and the quadratic fit methods when the same initial starting conditions are used for function \eqref{Rosenbrock2} with $\kappa=100$.

\begin{figure}[ht]
  \includegraphics[scale = 0.345]{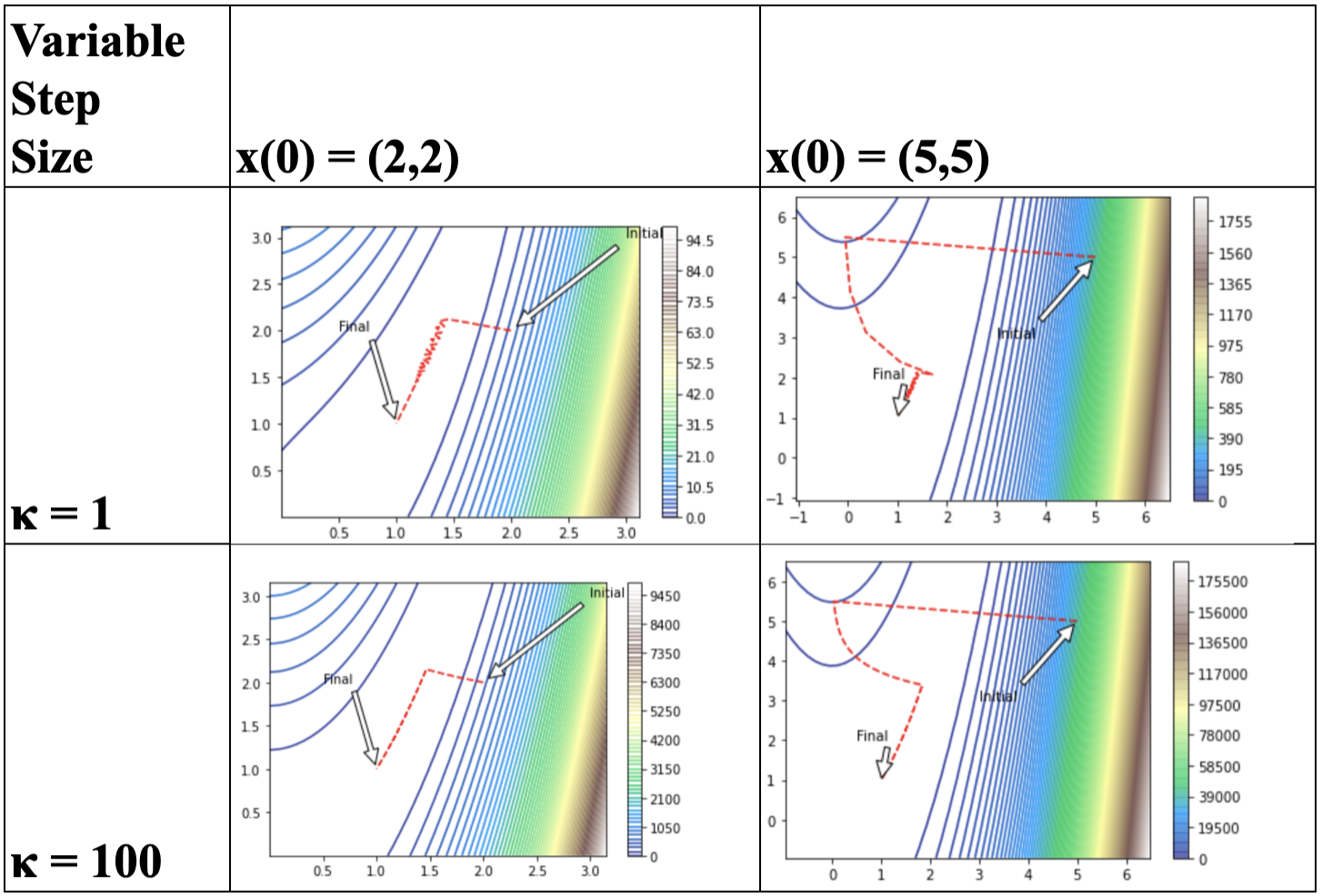}
  \caption{Level curves for Steepest Gradient Descent with variable step size.}
  \label{var}
\end{figure}

\begin{figure}[ht]
  \includegraphics[scale = 0.345]{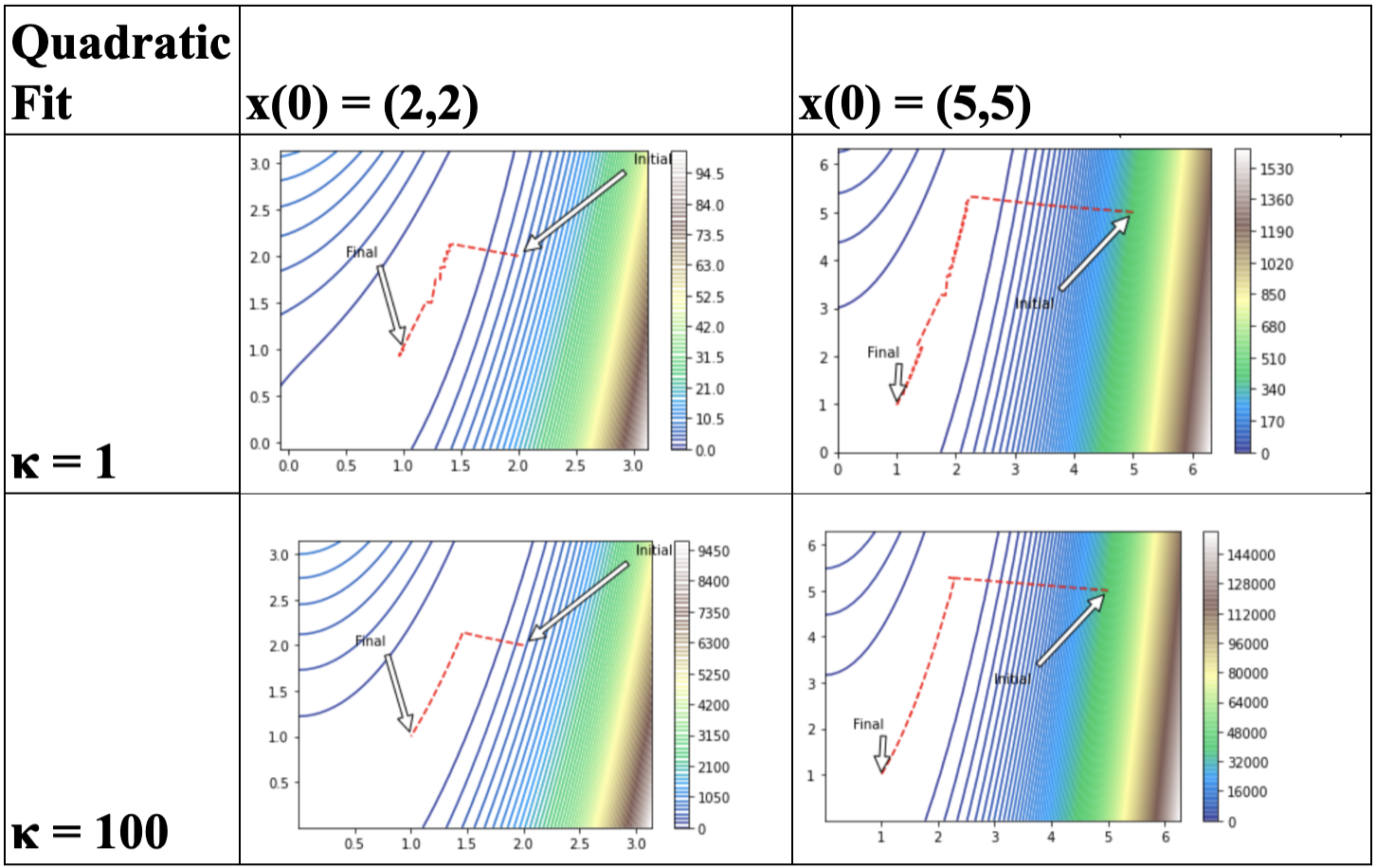}
  \caption{Level curves for Steepest Gradient Descent with quadratic fit.}
  \label{quad}
\end{figure}

\begin{figure}[ht]
  \includegraphics[scale = 0.345]{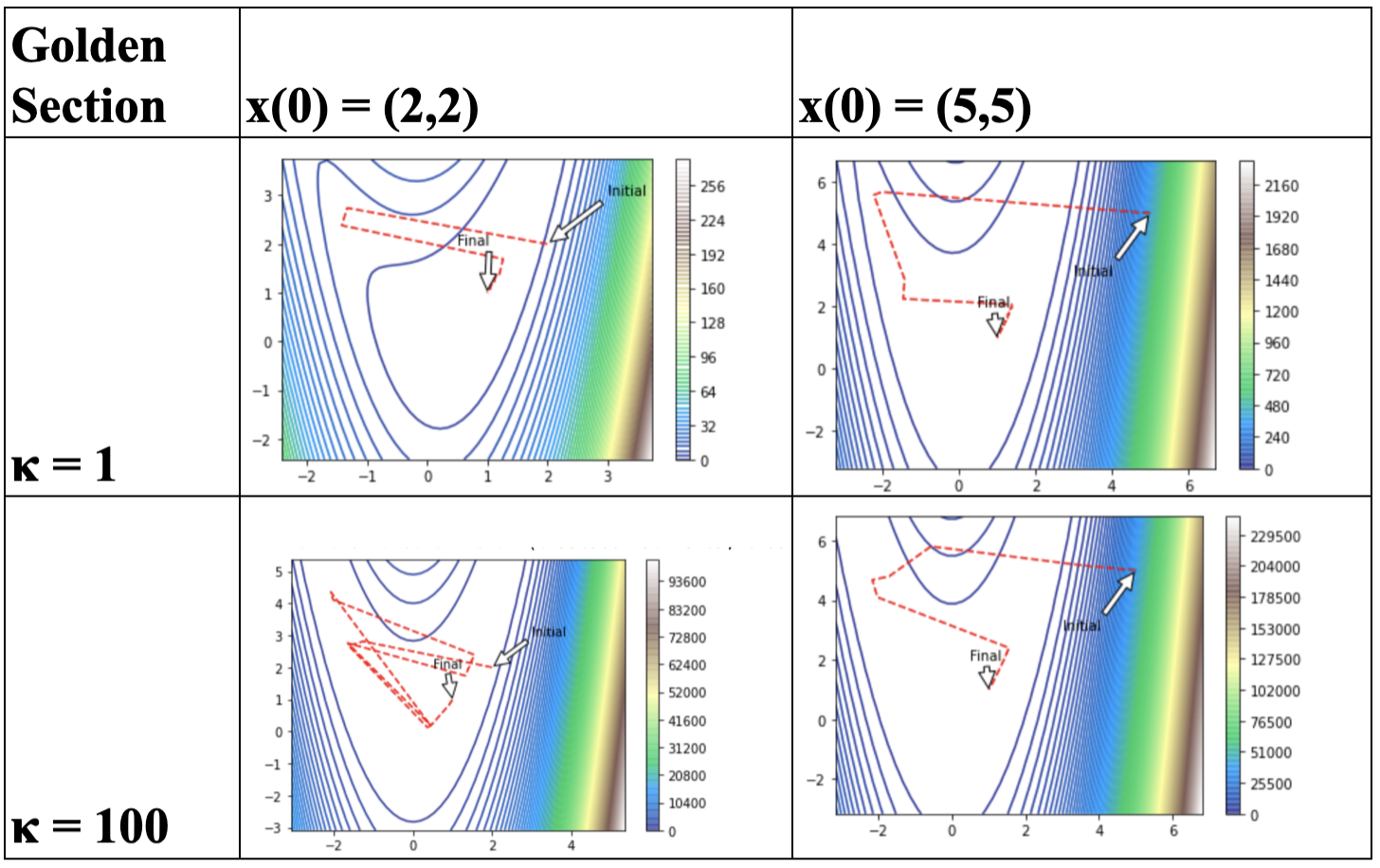}
  \caption{Level curves for Steepest Gradient Descent with golden section}
  \label{golden}
\end{figure}

\begin{figure}[ht]
  \includegraphics[scale = 0.262]{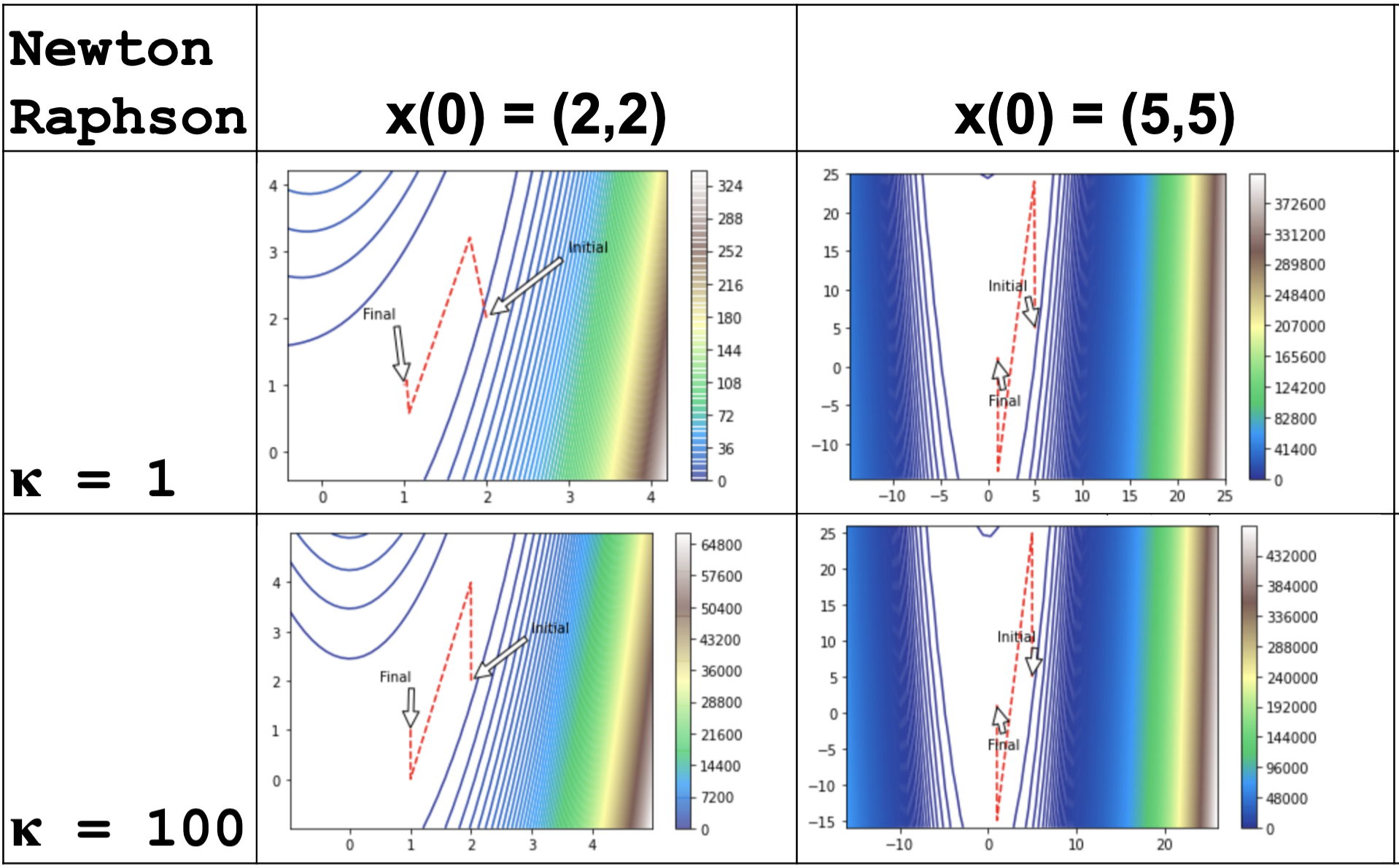}
  \caption{Level curves for Newton-Raphson's method.}
  \label{NR}
\end{figure}

\section{Conclusions}\label{sec:conclusion}
We analyze convergence attributes of some selected first and second order methods such as the steepest descent, Newton-Raphson, and conjugate gradient and apply it to a class of Rosenbrock functions.  We show through different minimization algorithms for function \eqref{Rosenbrock2} using values $\kappa=1$ and $\kappa=100$ that it is still possible for equation \eqref{Rosenbrock2} to converge to its minimum.Numerical experiments affirm that the Newton-Raphson method has the fastest convergence rate for the two strictly convex functions used in this paper provided the initial starting point is close to the minimum as seen with the starting points used.  To conclude, choosing the best method depends on the type of problem, the performance design specifications, and the resources available. As such, this study highlighted the differences and the trade-offs involved in comparing these algorithms to contribute to one's endeavor in selecting the most appropriate optimization method. 

\section*{ACKNOWLEDGEMENTS}
This work is done as part of a graduate course on Optimization Methods at the University of Central Florida.

\bibliographystyle{IEEEtran}
\small\bibliography{project.bib}




\end{document}